# A NOTE ON THE STATIONARY BOOTSTRAP'S VARIANCE

By Daniel J. Nordman

*Iowa State University*

Because the stationary bootstrap resamples data blocks of random length, this method has been thought to have the largest asymptotic variance among block bootstraps Lahiri [*Ann. Statist.* **27** (1999) 386–404]. It is shown here that the variance of the stationary bootstrap surprisingly matches that of a block bootstrap based on *nonrandom, nonoverlapping* blocks. This argument translates the variance expansion into the frequency domain and provides a unified way of determining variances for other block bootstraps. Some previous results on the stationary bootstrap, related to asymptotic relative efficiency and optimal block size, are also updated.

**1. Introduction.** Block bootstrap methods "recreate" an observed time series for inference by resampling and concatenating small blocks of data. Various block bootstraps exist and differ in their approaches to block resampling. Lahiri [9] investigated four block bootstrap methods for estimating the bias or variance of certain time series statistics based on sample means: namely, the moving block bootstrap (MBB) from Künsch [8] and Liu and Singh [11]; the nonoverlapping block bootstrap (NBB) from Carlstein [3]; the circular block bootstrap (CBB) from Politis and Romano [14]; and the stationary bootstrap (SB) from Politis and Romano [15]. The main findings in Lahiri [9] were that all four bootstrap methods have the same asymptotic bias but their asymptotic variances can differ due to the distinct blocking mechanism used in each method. For concreteness, we may consider a stationary, weakly dependent time series stretch $X_1, \ldots, X_n \in \mathbb{R}$ of size $n$ and estimate the variance $\sigma_n^2 \equiv n \operatorname{Var}(\bar{X}_n)$ of the sample mean $\bar{X}_n$ with a block bootstrap estimator. The resulting bootstrap estimator $\hat{\sigma}_\ell^2$ of $\sigma_n^2$ has a variance which behaves like

$$\operatorname{Var}(\hat{\sigma}_\ell^2) = C \frac{\ell}{n} + o(\ell/n) \qquad \text{as } n \to \infty, \tag{1}$$

where $\ell$ is the bootstrap block size used, $\ell^{-1} + \ell/n = o(1)$, and $C > 0$ is a constant that depends on the process spectral density as well as the block resampling device. Carlstein [3], Künsch [8], Bühlmann and Künsch [2] and Hall, Horowitz and Jing [7] had previously obtained some similar results for the MBB and the NBB. However, much development in Lahiri [9] was devoted to the SB in particular, because this method resamples blocks of *random* length and thereby differs from the other block bootstraps above using fixed block lengths $\ell$. Specifically, the variance expansions in [9] (Theorem 3.2) suggest that the constants $C$ in (1) differ across block bootstrap estimators so that

$$
\lim_{n \to \infty} \frac{\mathrm{Var}(\hat{\sigma}^2_{\ell,CBB})}{\mathrm{Var}(\hat{\sigma}^2_{\ell,NBB})} = \lim_{n \to \infty} \frac{\mathrm{Var}(\hat{\sigma}^2_{\ell,MBB})}{\mathrm{Var}(\hat{\sigma}^2_{\ell,NBB})} = \frac{2}{3} \quad \text{and}
$$
(2)
$$
\lim_{n \to \infty} \frac{\mathrm{Var}(\hat{\sigma}^2_{\ell,SB})}{\mathrm{Var}(\hat{\sigma}^2_{\ell,NBB})} > 2.
$$

We claim however that, contrary to (2), the SB actually has the *same* large sample variance as the NBB method based on nonrandom, nonoverlapping blocks (to the first order). For illustration, Figure 1 displays the ratio of SB and NBB estimator variances $\mathrm{Var}(\hat{\sigma}^2_{\ell,SB})/\mathrm{Var}(\hat{\sigma}^2_{\ell,NBB})$, approximated by simulation, with samples from three different AR(1) models. As the sample sizes increase, the ratio of variances tends to 1 under each process model. This result is quite unexpected because the SB estimator had been anticipated to have the largest variance among the block bootstraps due to "extra randomness" induced by resampling blocks of random length.

Additionally, Politis and White [16] raised concerns on past claims regarding the potential size of the SB variance and the asymptotic relative efficiency (ARE) of the SB to other block bootstraps. This manuscript corrects the variance of the SB estimator, and thereby invalidates all prior comparisons of ARE involving the SB as well as previous optimal block size calculations for SB (cf. [9, 10, 16]). We shall update these results in Section 2.2 to follow and modify a proposal of Politis and White [16] for accurately estimating optimal SB block sizes.

In reformulating the variance of the SB estimator, we also provide a new proof of this result. To determine the leading component in the SB variance, the proof in Lahiri [9] involves conditioning arguments and moment bounds for random block sums as well as Hilbert space properties. The approach given here transfers the SB variance expansion into the frequency domain and uses cumulant expressions for the periodogram [1, 5]. This reduces the complexity of the arguments and more readily isolates the main contributor to the SB variance. The frequency domain arguments appear to be new and also apply for determining variances of other block bootstrap estimators in a simple manner.

We end this section by briefly mentioning other block bootstrap methods. Section 2 describes the SB variance estimator for the sample mean and provides its large sample variance. Section 3 gives a general variance result for lag-weighted sums of sample covariances, which can be used to determine the variances of the SB and other block bootstraps. Section 4 proves this result.

Additional block bootstrap methods include the matched block bootstrap due to Carlstein, Do, Hall, Hesterberg and Künsch [4] and the tapered block bootstrap (TBB) proposed by Paparoditis and Politis [12]. These have similar order variances but smaller biases than the other block bootstraps mentioned here, though the matched block bootstrap requires stronger assumptions on the time process (e.g., Markov) than the TBB for bias reduction. See [9, 10] for further references on the block bootstrap and block estimation.

## 2. Large sample variance of the stationary bootstrap.

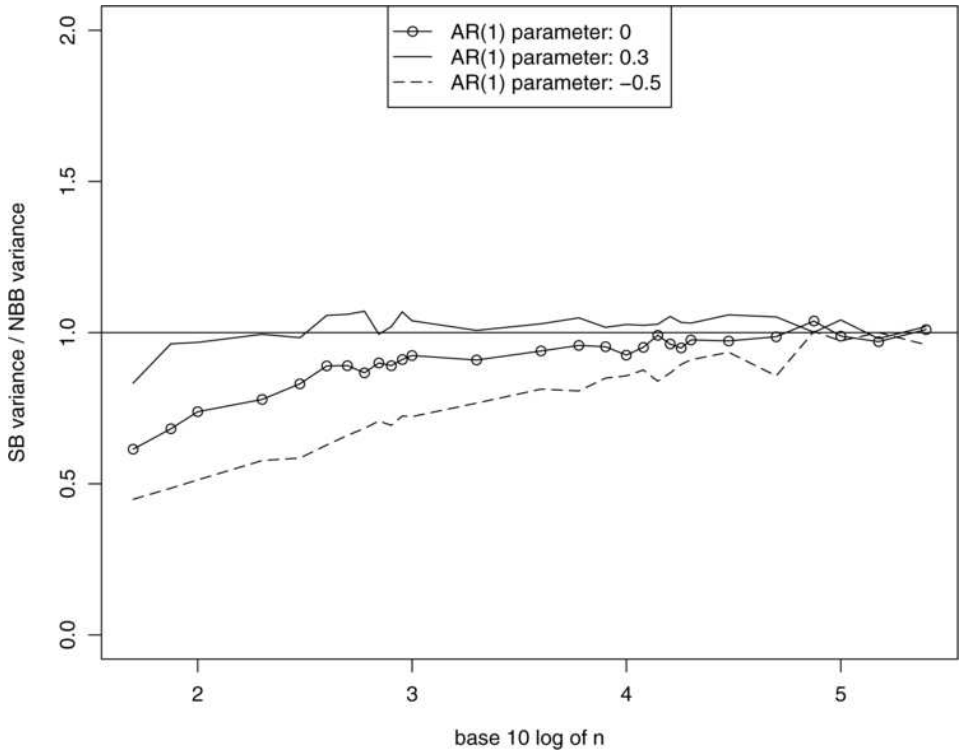

FIG. 1. *Ratio of SB and NBB variances,* $\text{Var}(\hat{\sigma}^2_{\ell,SB})/\text{Var}(\hat{\sigma}^2_{\ell,NBB})$, *for block* $\ell = n^{1/3}$ *and sample sizes n from AR(1) models:* $X_t = \phi X_{t-1} + Z_t$, $\{Z_t\}$ *i.i.d. standard normals,* $\phi = 0, 0.3, -0.5$ *(variances based on 5000 simulations). See Section 2.1 for these estimators.*

2.1. *Stationary bootstrap estimator.* We suppose that the data $X_1, \ldots, X_n$ arise from a real-valued strictly stationary process $\{X_t\}$, $t \in \mathbb{Z}$. For simplicity, we shall focus the exposition on the bootstrap estimator for the variance $\sigma_n^2 = n\operatorname{Var}(\bar{X}_n)$ of the sample mean $\bar{X}_n$, but similar results hold for variance estimation of "nearly linear" statistics formed by a smooth function of sample means; see Remark 1 in Section 2.2.

The SB of Politis and Romano [15] generates a bootstrap version of the observed data $X_1, \ldots, X_n$ by resampling data blocks of random length. To describe the method, denote a data block as $B(i, k) = (X_i, \ldots, X_{i+k-1})$ for $i, k \geq 1$ where the data are periodically extended as $X_i = X_{i-n}$ for $i > n$. Let $J_1, \ldots, J_n$ be i.i.d. geometric random variables, independent of the data, where the event "$J_1 = k$," $k \geq 1$, has probability $pq^{k-1}$ for $p \equiv \ell^{-1} \in (0, 1)$, $q = 1 - p$. For $K = \min\{k : \sum_{i=1}^k J_i \geq n\}$, let $I_1, \ldots, I_K$ be i.i.d. uniform variables on $\{1, \ldots, n\}$. Then, the SB sample $X_1^*, \ldots, X_n^*$ is given by the first $n$ observations in the sequence $B(I_1, J_1), \ldots, B(I_K, J_K)$ and $\ell = p^{-1}$ represents the expected length of a resampled block. Under this resampling, Politis and Romano [15] show that the SB sample exhibits stationarity. The SB estimator of $\sigma_n^2$ is the variance $\hat{\sigma}_{\ell,SB}^2 \equiv n\operatorname{Var}_*(\bar{X}_n^* | X_1, \ldots, X_n)$ of the bootstrap average $\bar{X}_n^* = n^{-1}\sum_{i=1}^n X_i^*$ conditional on the data, which has a closed-form expression in this case ([15], Lemma 1) as

$$(3) \quad \hat{\sigma}_{\ell,SB}^2 = \hat{r}(0) + 2\sum_{k=1}^{n-1} q_{k,n} \hat{r}(k), \qquad q_{k,n} = (1 - n^{-1}k)q^k + n^{-1}k q^{n-k}$$

based on sample covariances $\hat{r}(k) = n^{-1}\sum_{t=1}^{n-k}(X_t - \bar{X}_n)(X_{t+k} - \bar{X}_n)$.

For comparison, we recall the NBB variance estimator which resamples disjoint blocks of fixed integer length $\ell \geq 1$ as suggested by Carlstein [3]. For this, we randomly sample $b = \lfloor n/l \rfloor$ blocks with replacement from the nonoverlapping collection $\{B(\ell(i-1)+1, \ell i) : i = 1, \ldots, b\}$ and concatenate these to produce the NBB sample $X_{1,NBB}^*, \ldots, X_{b\ell,NBB}^*$ of size $b\ell$ with sample mean $\bar{X}_{n,NBB}^* = (b\ell)^{-1}\sum_{i=1}^{b\ell} X_{i,NBB}^*$. The NBB variance estimator for $\sigma_n^2$ is

$$\hat{\sigma}_{\ell,NBB}^2 \equiv bl\operatorname{Var}_*(\bar{X}_{n,NBB}^* | X_1, \ldots, X_n) = \frac{\ell}{b}\sum_{i=1}^b (\bar{X}_{i,NBB} - \tilde{X}_n)^2,$$

which is a scaled sample variance of block means $\bar{X}_{i,NBB} \equiv \ell^{-1}\sum_{k=\ell(i-1)+1}^{\ell i} X_k$, $i = 1, \ldots, b$ with $\tilde{X}_n \equiv b^{-1}\sum_{i=1}^b \bar{X}_{i,NBB}$.

2.2. *Variance expression for stationary bootstrap estimator.* We state the SB variance in Theorem 1, which requires some notation. Denote the covariances and spectral density of $\{X_t\}$ as $r(k) = \operatorname{Cov}(X_0, X_k)$, $k \in \mathbb{Z}$ and $f(\omega) = (2\pi)^{-1} \sum_{k \in \mathbb{Z}} r(k) e^{-\imath k \omega}$, $\omega \in \mathbb{R}$, $\imath = \sqrt{-1}$. Let $\mathrm{cu}(Y_1, \ldots, Y_m)$ denote the cumulant of arbitrary random variables $(Y_1, \ldots, Y_m)$, $m \geq 1$. The assumptions are A.1–A.3 as follows:

A.1: $|r(0)| + \sum_{k=1}^{\infty} k|r(k)| < \infty;$ A.2: $\sum_{t_1, t_2, t_3 \in \mathbb{Z}} |\mathrm{cu}(X_0, X_{t_1}, X_{t_2}, X_{t_3})| < \infty;$

and A.3: $\ell \to \infty$ and $\ell \log n / n \to 0$ as $n \to \infty$. Assumptions A.1–A.2 essentially correspond to those of Politis and Romano [15] for showing consistency of the SB estimator $\hat{\sigma}^2_{\ell, SB}$. These hold under the mixing-based moment conditions assumed by Künsch [8], Politis and White [16] and Paparoditis and Politis [12] for stating variance expansions for MBB/NBB, CBB and TBB estimators, respectively, of the sample mean's variance. The SB variance result will be derived in Section 3, requiring a slightly stronger block condition A.3 than $\ell/n = o(1)$ introduced by Künsch [8] (Theorem 3.3) for MBB variance calculations. Lemma 1 of Section 4 describes how A.1 is applied.

For completeness, we include the bias of the SB estimator in Theorem 1, originally given by [15], in (9).

THEOREM 1. *Under assumptions* A.1–A.3 *and as* $n \to \infty$, *the variance of the SB estimator* $\hat{\sigma}^2_{\ell, SB}$ *of* $\sigma^2_n = n \mathrm{Var}(\bar{X}_n)$ *is*

$$\mathrm{Var}(\hat{\sigma}^2_{\ell, SB}) = \frac{\ell}{n} 2 \{2\pi f(0)\}^2 + o(\ell/n).$$

*Under* A.1 *and* $\ell^{-1} + \ell^2/n = o(1)$ *as* $n \to \infty$, *the bias of the SB estimator is*

$$\mathrm{E}(\hat{\sigma}^2_{\ell, SB}) - \sigma^2_n = -\frac{G}{\ell} + o(1/\ell) \qquad where \ G = \sum_{k \in \mathbb{Z}} |k| r(k).$$

The SB variance and bias expansions in Theorem 1 match those of the NBB estimator for $\sigma^2_n$ (cf. [3, 9]). Hence, while both SB and NBB variance estimators appear structurally very different, these estimators have the same asymptotic performance in terms of mean squared error (MSE) as

(4) $\quad MSE(\hat{\sigma}^2_{\ell, SB}) = \frac{\ell}{n} 2\{2\pi f(0)\}^2 + \frac{G^2}{\ell^2} + o(\ell/n) + o(\ell^{-2})$

matches the expression for $MSE(\hat{\sigma}^2_{\ell, NBB})$. Theorem 1 also helps to explain some simulation evidence in Politis and White [16], which indicated that the finite sample MSE performance of the SB estimator was better than predicted by previous asymptotic theory. We mention that, as noted in [4], the matched block bootstrap estimator of $\sigma^2_n$ has the same first-order variance as the NBB as well. Interestingly, the matched block bootstrap seems at first glance to introduce more variability through its block resampling scheme (i.e., block "matching") but this aspect does not inflate this bootstrap estimator's variance as one might expect, much like the SB.

The next section describes the main component in the proof of Theorem 1.

REMARK 1. Although we focus on SB variance estimation for a real-valued sample mean, Theorem 1 also applies for more complicated statistics formed by a smooth function $H(\bar{Y}_n)$ of the sample mean $\bar{Y}_n$ of a multivariate time series $Y_1, \ldots, Y_n \in \mathbb{R}^d$, $H: \mathbb{R}^d \to \mathbb{R}$ (cf. [6, 9]). The variance of the SB estimator for $n \operatorname{Var}(H(\bar{Y}_n))$ is typically determined by variance of the SB estimator for $\sigma_n^2 = n \operatorname{Var}(\bar{X}_n)$ based on the linearly approximating series $X_1, \ldots, X_n$, with $X_t = H(\mu) + c'(Y_t - \mu) \in \mathbb{R}$ involving the first partial derivatives $c$ of $H$ at $\mu \equiv \mathrm{E}(Y_t)$ (cf. [7, 9]).

REMARK 2. We discuss the apparent difference between SB variance results in Lahiri [9] and those presented here. The SB variance expression in [9] was determined by the limit of a quantity $\lim_{n\to\infty} \int_{-\pi}^{\pi} e^{\iota\omega} g(\omega)/(1 - qe^{\iota\omega}) \, d\omega$ [see Lemma 5.3(i) there], where $g(\omega) = 4f^2(\omega)$ for the case of estimating $\sigma_n^2 = n \operatorname{Var}(\bar{X}_n)$. In the proof of this limit, "$-(\sin\omega)^2$" should seemingly replace "$(\sin\omega)^2$" in lines (5), (8) of [9], which produces $\lim_{n\to\infty} \int_{-\pi}^{\pi} e^{\iota\omega} g(\omega)/(1 - qe^{\iota\omega}) \, d\omega = \pi g(0) - 2^{-1} \int_{-\pi}^{\pi} g(\omega) \, d\omega$. With this change, the proof in [9] gives results matching those here for the sample mean; further details can be found in [10] (pp. 138, 144). More generally, the results and technical arguments in [9] remain valid for block bootstrap estimation of bias and variance by changing the variances of the SB to match those provided for the NBB.

REMARK 3. Politis and Romano [15] showed that an expected block size of form $\ell_{SB}^{\mathrm{opt}} = Cn^{1/3}$ asymptotically minimizes the MSE (4) of the SB variance estimator and Theorem 1 yields the optimal constant as $C = |G/\{2\pi f(0)\}|^{2/3}$ assuming $G \neq 0$, $f(0) > 0$. This block expression also replaces the one in (6) of Politis and White [16] and their Section 3.2 method of estimating $\ell_{SB}^{\mathrm{opt}}$ requires slight modification; their procedure would now produce a block estimator $\hat{\ell}_{SB}^{\mathrm{opt}} \equiv (2/3)^{1/3} \hat{\ell}_{CBB}^{\mathrm{opt}}$ using their Section 3.3 estimator $\hat{\ell}_{CBB}^{\mathrm{opt}}$ of the optimal CBB block length (which follows from a relationship between these optimal blocks described in Remark 4). These changes in Politis and White [16] follow by correcting their Theorem 3.1 (based on Lahiri [9]) so that "$D_{SB} \equiv (3/2) D_{CBB}$" there, which reproduces the SB variance as in Theorem 1.

REMARK 4. Given the SB variance, we may update a result (Lemma 3.1) from Politis and White [16] on the asymptotic relative efficiency (ARE) of SB versus CBB methods in estimating $\sigma_n^2 = n \operatorname{Var}(\bar{X}_n)$. In discussing ARE, we may also substitute "NBB" for "SB" or "MBB" for "CBB" because the estimators in this pairing match in variance and optimal block size to the first order; all four estimators have the same bias expression as given in Theorem 1. Compared to the SB, the CBB estimator has a smaller asymptotic variance by a factor of $2/3$ (see Section 3) which produces a longer MSE-optimal block length as $\ell^{\mathrm{opt}}_{CBB} = (3/2)^{1/3}\ell^{\mathrm{opt}}_{SB}$ by modifying (4). At each method's respective optimal block size, the ARE of the SB to the CBB estimator is given by the limiting ratio of minimized MSEs

$$ARE(\hat{\sigma}^2_{\ell^{\mathrm{opt}},SB}; \hat{\sigma}^2_{\ell^{\mathrm{opt}},CBB}) = \lim_{n\to\infty} \frac{MSE(\hat{\sigma}^2_{\ell^{\mathrm{opt}},CBB})}{MSE(\hat{\sigma}^2_{\ell^{\mathrm{opt}},SB})} = (2/3)^{2/3},$$

which replaces the previous bounds on this ARE given in Lemma 3.1 of [16].

**3. A general variance result.** The SB variance in Theorem 1 follows from a general variance result, provided in Theorem 2 below, for weighted sums $T_n$ of sample covariances $\hat{r}(k)$ given by

$$(5) \qquad T_n \equiv \sum_{k=0}^{n-1} a_{k,n}\hat{r}(k) \qquad \text{with } A_n \equiv \sum_{k=1}^{n-1} a_{k,n}^2(1-n^{-1}k)^2 > 0,$$

based on a real-valued stationary sample $X_1,\ldots,X_n$ and a triangular array of weights $\{a_{0,n},\ldots,a_{n-1,n}\}\subset\mathbb{R}$, $n\geq 1$. While several authors have considered large sample variances of certain lag weight spectral estimates (cf. Theorem 5.6.2/5.9.11 and page 161 references, [1]; Theorem 5A, [13]), these results do not readily fit the framework (3) required for the SB variance.

Define a smoothing window as $H_n(\omega) \equiv \sum_{k=1}^{n-1} a_{k,n}(1-n^{-1}k)e^{-\imath k\omega}$ and a nonnegative kernel $K_n(\omega) = H_n(\omega)H_n(-\omega)/(2\pi A_n)$, $\omega\in\mathbb{R}$. Let

$$B_n \equiv n^{-1} + \left(\sum_{k=1}^{n-1}|a_{k,n}|\right)^2 n^{-2}\log n + n^{-1}\sum_{k=1}^{n-1}|a_{k,n}|n^{-1}k(1-n^{-1}k).$$

THEOREM 2. *Suppose A.1–A.2 hold and $\sup_n \max_{0\leq k\leq n-1}|a_{k,n}| < \infty$. Then:*

(i) $\mathrm{Var}(T_n) = \dfrac{A_n(2\pi)^2}{n}\displaystyle\int_{-\pi}^{\pi} K_n(\omega)f^2(\omega)\,d\omega + O((A_nB_n/n)^{1/2} + B_n),$

(ii) $\displaystyle\lim_{n\to\infty}\int_{-\pi}^{\pi} K_n(\omega)f^2(\omega)\,d\omega = f^2(0)$

*as $n\to\infty$, assuming additionally for* (ii) *that $\lim_{n\to\infty} A_n = \infty$ and $\sup_n \sum_{k=2}^{n-1}|a_{k,n} - a_{k-1,n}| < \infty$.*

Section 4 provides a proof of Theorem 2. The factor "$(1-n^{-1}k)$" is important in the definition of $A_n$ from (5) and nullifies weights $a_{k,n}$ at large lags $k$ in the sum. This becomes relevant in the computation of the SB variance as well as the CBB variance, as discussed below.

The SB variance estimator (3) has a sum form (5) with weights $a_{0,n} = 1$ and $a_{k,n} = 2q_{k,n} = 2\{(1-n^{-1}k)q^k + n^{-1}kq^{n-k}\}$ for $k \geq 1$ so that

$$A_n = 4 \sum_{k=1}^{n-1} q^{2k} + O(\ell^2/n) = \frac{4q^2}{1-q^2} + O(\ell^2/n), \qquad B_n = O(\ell^2 \log n/n^2 + n^{-1}),$$

recalling $1 - q = p = \ell^{-1}$ and using $\sum_{k=1}^{n-1} q_{k,n} \leq 2p^{-1}$, $\sum_{k=1}^{\infty} kq^{2k} \leq p^{-2}$ and $\log q \leq -p$. The terms "$n^{-1}kq^{n-k}$" in $q_{k,n}$ can essentially be ignored in computing $A_n$ because of the factor "$(1-n^{-1}k)$" in (5). The SB variance in Theorem 1 follows immediately from Theorem 2 with the block assumption A.3 since $\sum_{k=2}^{n-1} |q_{k,n} - q_{k-1,n}| \leq 2$, $\lim_{n \to \infty} \ell^{-1} A_n = 2$ and $B_n = o(\ell/n)$.

Theorem 2 also applies for determining the variance of other block bootstrap estimators that use overlapping blocks. For example, CBB, MBB and TBB estimators of $\sigma_n^2 = n \operatorname{Var}(\bar{X}_n)$, based on a fixed block length $\ell$, can also be written as weighted sums $T_n$ of covariance estimators with weights $a_{0,n} = 1$ and

CBB: $a_{k,n} = 2(1 - \ell^{-1}k)$ for $1 \leq k < \ell$, $\qquad a_{k,n} = 0$ for $\ell \leq k \leq n - \ell$,

$\qquad a_{k,n} = a_{n-k,k}$ for $n - \ell < k < n$;

MBB: $a_{k,n} = 2(1 - \ell^{-1}k)$ for $1 \leq k < \ell$, $\qquad a_{k,n} = 0$ for $k \geq \ell$;

TBB: $a_{k,n} = 2v_n(k)/v_n(0)$ for $1 \leq k < \ell$, $\qquad a_{k,n} = 0$ for $k \geq \ell$,

where TBB weights $v_n(k)$ appear in (3.6) of Künsch [8] (depending on a tapering window of bounded variation). The fact that the TBB estimator can be expressed as a sum (5) follows by the equivalence of the TBB and the tapered block jackknife estimator of $\sigma_n^2$ (Lemma A.1, [12]) combined with Theorem 3.9 of [8] which expresses this jackknife estimator as a lag-weighted sum of sample covariances. The MBB is a special, untapered case of the TBB where $v_n(k) = \ell - k$, $k < \ell$. Technically, the MBB/TBB do not precisely use sample covariances $\hat{r}(k)$ as defined here (see (3.7), (3.8) of [8]) but their covariance estimators are close so that any differences have a negligible $o(\ell/n)$ impact on the variance expansions of the MBB/TBB estimators. The CBB estimator $\hat{\sigma}_{\ell,CBB}^2$ of $\sigma_n^2$ has a sum form (5) because

$$\hat{\sigma}_{\ell,CBB}^2 = \hat{r}(0) + 2 \sum_{k=1}^{\ell-1} (1 - \ell^{-1}k)\hat{r}_{\mathrm{cir}}(k), \qquad \hat{r}_{\mathrm{cir}}(k) = \hat{r}(k) + \hat{r}(n-k)$$

using circular sample covariances $\hat{r}_{\mathrm{cir}}(k) \equiv n^{-1} \sum_{t=1}^{n}(X_t - \bar{X}_n)(X_{t+k} - \bar{X}_n)$.

Considering the MBB, for example, it is easy to check that $\ell^{-1} A_n = 4\ell^{-1} \times \sum_{k=1}^{\ell}(1-\ell^{-1}k)^2 + O(\ell/n) \to 4/3$ and $B_n = o(l/n)$ as $n \to \infty$ under the block condition A.3, in which case Theorem 2 produces the well-known variance of the MBB estimator $\text{Var}(\hat{\sigma}^2_{\ell,MBB}) = (4/3)\{2\pi f(0)\}^2 \ell/n + o(\ell/n)$ (cf. [7, 8]). For the CBB estimator $\hat{\sigma}^2_{\ell,CBB}$, this same variance expression and $\ell^{-1} A_n \to 4/3$ both hold and, in this case, the factor "$(1 - n^{-1}k)$" in the definition of $A_n$ wipes out the effect of CBB weights $a_{k,n}$, $k > n - \ell$. The variance of the TBB estimator can be similarly replicated with Theorem 2; see Theorem 3.3 of [8] or Theorem 2 of [12] for this expression.

**4. Proof of Theorem 2.** We require some notation in the frequency domain. Let $I_n(\omega) = d_n(\omega) d_n(-\omega)/\{2\pi n\}$ denote the periodogram based on the discrete Fourier transform $d_n(\omega) = \sum_{t=1}^n X_t e^{-\imath t \omega}$, $\omega \in \mathbb{R}$. Define the Fejer kernel $K_n^F(\omega) = H_n^F(\omega) H_n^F(-\omega)/\{2\pi n\} = (2\pi)^{-1} \sum_{t=-n}^n (1 - n^{-1}|t|) e^{-\imath t \omega}$ using the spectral window $H_n^F(\omega) = \sum_{t=1}^n e^{-\imath t \omega}, \omega \in \mathbb{R}$. Define also a function from Dahlhaus [5] given by

$$L_{0,n}(\omega) \equiv \begin{cases} n, & |\omega| \leq 1/n, \\ 1/|\omega|, & 1/n < |\omega| \leq \pi, \end{cases}$$

with periodic extension $L_{0,n}(\omega) = L_{0,n}(2\pi + \omega)$, $\omega \in \mathbb{R}$, which is helpful for bounding both $H_n^F(\cdot)$ and the spectral window $H_n(\cdot)$ from Theorem 2. Throughout the remainder, let $C > 0$ denote a generic constant which does not depend on $n$ or any frequencies $\omega, \omega_1, \omega_2, \lambda \in \mathbb{R}$. We summarize some technical facts in the following Lemma 1. Regarding Lemma 1(ii), we add that uniform bounds on $|\text{cu}(d_n(\omega_1), d_n(\omega_2)) - 2\pi H_n^F(\omega_1 + \omega_2) f(\omega_1)|$ may range from $o(n)$ when $\sum_{k \in \mathbb{Z}} |r(k)| < \infty$ to $O(1)$ under A.1, with intermediate possibilities; see Theorem 1 of [5]. Similarly to Theorem 5.6.2 of [1], we require A.1 for tractability due to weak conditions on the weight forms $a_{k,n}$ in (5).

LEMMA 1. *Let* $\omega_1, \omega_2 \in \mathbb{R}$. *Under assumption* A.1 *for parts* (ii)–(iii) *below:*

(i) $|H_n^F(\omega_1)| \leq C L_{0,n}(\omega_1)$ *and* $\int_{-\pi}^{\pi} L_{0,n}(\omega) d\omega \leq C \log n$.
(ii) $\text{cu}(d_n(\omega_1), d_n(\omega_2)) = 2\pi H_n^F(\omega_1 + \omega_2) f(\omega_1) + R_n(\omega_1, \omega_2)$, $|R_n(\omega_1, \omega_2)| \leq C$.
(iii) $\sum_{d=1}^{\infty} d \sum_{t \in \mathbb{Z}} |r(t) r(t+d)| \leq C$.

PROOF. Part (i) is from [5], Lemma 1. Part (ii) follows from Theorem 4.3.2 of [1] under A.1. Part (iii) follows from A.1 and $(1 + |t|)(1 + |t + d|) > d$ for $d \geq 1$, $t \in \mathbb{Z}$. □

Because sample covariances $\hat{r}(k)$ are defined with a sample mean correction as in (3), we may assume that $\mathbb{E} X_t = 0$ holds without loss of generality. Then, sample covariances defined with a process mean $\mathbb{E} X_t$ correction are $\tilde{r}(k) = n^{-1} \sum_{t=1}^{n-k} X_t X_{t+k}$, $0 \leq k < n$ so that $\hat{r}(k) = \tilde{r}(k) - n^{-1} \bar{X}_n \sum_{t=k+1}^{n-k} X_t - n^{-1} k \bar{X}_n^2$. We may decompose (5) as

$$T_n = T_{1n} + a_{0,n} \hat{r}(0) - \tilde{T}_{1n} - \tilde{T}_{2n}, \qquad T_{1n} \equiv \sum_{k=1}^{n-1} a_{k,n} \tilde{r}(k),$$

$\tilde{T}_{1n} \equiv n^{-1} \bar{X}_n \sum_{k=1}^{n-1} a_{k,n} \sum_{t=k+1}^{n-k} X_t$ and $\tilde{T}_{2n} \equiv n^{-1} \bar{X}_n^2 \sum_{k=1}^{n-1} k a_{k,n}$. Note that $\mathrm{Var}(T_n - T_{1n}) = O(B_n)$ since $\mathrm{Var}(\hat{r}(0)) = O(n^{-1})$ and $\mathrm{Var}(\tilde{T}_{in}) = O(\{n^{-1} \times \sum_{k=1}^{n-1} |a_{k,n}|\}^2)$ for $i=1,2$ under A.1, A.2. We aim to show

$$(6) \quad \mathrm{Var}(T_{1n}) = M_n + O(B_n), \qquad M_n \equiv \frac{A_n(2\pi)^2}{n} \int_{-\pi}^{\pi} K_n(\omega) f^2(\lambda) \, d\omega$$

so that Theorem 2(i) will then follow by the Cauchy–Schwarz inequality since $f(\cdot)$ is bounded by A.1 and $K_n(\cdot)$ is nonnegative with $\int_{-\pi}^{\pi} K_n(\omega) \, d\omega = 1$.

We transfer to the frequency domain using that $\tilde{r}(k) = \int_{-\pi}^{\pi} e^{-\imath k\omega} I_n(\omega) \, d\omega$, $1 \leq k \leq n-1$ to rewrite $T_{1n} = \sum_{k=1}^{n-1} a_{k,n} \int_{-\pi}^{\pi} e^{-\imath k\omega} I_n(\omega) \, d\omega$ and find

$$\mathrm{Var}(T_{1n}) = \sum_{j,k=1}^{n-1} a_{j,n} a_{k,n} \int_{-\pi}^{\pi} \int_{-\pi}^{\pi} e^{-\imath j\omega} e^{\imath k\lambda} \mathrm{cu}\{I_n(\omega), I_n(\lambda)\} \, d\lambda \, d\omega$$

$$(7)$$

$$= \sum_{i=1}^{3} S_{in}, S_{in} \equiv \sum_{j,k=1}^{n-1} \frac{a_{j,n} a_{k,n}}{(2\pi n)^2} \int_{-\pi}^{\pi} \int_{-\pi}^{\pi} e^{-\imath j\omega} e^{\imath k\lambda} g_i(\omega, \lambda) \, d\lambda \, d\omega,$$

since the cumulant product theorem (Theorem 2.3.2, [1]) and $\mathrm{E} X_t = 0$ entail

$$(2\pi n)^2 \mathrm{cu}\{I_n(\omega), I_n(\lambda)\} = \sum_{i=1}^{3} g_i(\omega, \lambda), \qquad \omega, \lambda \in \mathbb{R},$$

for $g_1(\omega, \lambda) \equiv \mathrm{cu}\{d_n(\omega), d_n(-\lambda)\} \mathrm{cu}\{d_n(-\omega), d_n(\lambda)\}$, $g_2(\omega, \lambda) \equiv g_1(\omega, -\lambda)$ and $g_3(\omega, \lambda) \equiv \mathrm{cu}\{d_n(\omega), d_n(-\omega), d_n(\lambda), d_n(-\lambda)\}$.

By A.2 and $\sup_{k,n} |a_{k,n}| \leq C$, it follows that $|S_{3n}| = O(1/n)$. We then use Lemma 1 to simplify the expressions for $S_{1n}$ and $S_{2n}$ in (7). In particular, by applying Lemma 1(ii) with the definitions of the kernel $K_n^F(\cdot)$ and smoothing window $H_n^F(\cdot)$, we may write $S_{1n} = \tilde{S}_{1n} + R_{1n}$, where

$$\tilde{S}_{1n} = \frac{2\pi}{n} \sum_{j,k=1}^{n-1} a_{j,n} a_{k,n} \int_{-\pi}^{\pi} \int_{-\pi}^{\pi} e^{-\imath j\omega} e^{\imath k\lambda} K_n^F(\omega - \lambda) f^2(\omega) \, d\lambda \, d\omega$$

$$(8)$$

$$= \frac{2\pi}{n} \sum_{j,k=1}^{n-1} a_{j,n} a_{k,n} (1 - n^{-1} k) \int_{-\pi}^{\pi} e^{-\imath (j-k)\omega} f^2(\omega) \, d\omega$$

using $\int_{-\pi}^{\pi} e^{\imath k\lambda} K_n^F(\omega - \lambda)\, d\lambda = (1 - n^{-1}k) e^{\imath k\omega}$, and we apply Lemma 1(i)–(ii) to bound the remainder

$$|R_{1n}| \leq \frac{C}{n^2}\left(\sum_{k=1}^{n-1} |a_{k,n}|\right)^2 \int_{-\pi}^{\pi} \int_{-\pi}^{\pi} \{L_{0,n}(\omega - \lambda) + 1\}\, d\lambda\, d\omega \leq CB_n.$$

A similar expansion $S_{2n} = \tilde{S}_{2n} + R_{2n}$ holds where $R_{2n} = O(B_n)$ and $\tilde{S}_{2n}$ is defined by replacing "$(j-k)$" with "$(j+k)$" in the exponent of (8). Since

$$(9) \qquad \int_{-\pi}^{\pi} e^{-\imath d\omega} f^2(\omega)\, d\omega = \frac{1}{2\pi} \sum_{t \in \mathbb{Z}} r(t) r(t+d), \qquad d \in \mathbb{Z},$$

we use a substitution $j + k = d$ in the sum $\tilde{S}_{2n}$ and arrange terms to find

$$|\tilde{S}_{2n}| \leq \frac{C}{n} \sum_{d=2}^{2n-2} d \sum_{t \in \mathbb{Z}} |r(t) r(t+d)| \leq \frac{C}{n}$$

by Lemma 1(iii). Finally, from (8) we have $\tilde{S}_{1n} = M_n + R_{3n}$ for $M_n$ in (6) and a remainder

$$R_{3n} \equiv \frac{2\pi}{n} \sum_{j,k=1}^{n-1} a_{j,n} a_{k,n} (1 - n^{-1}k) n^{-1} j \int_{-\pi}^{\pi} e^{-\imath(j-k)\omega} f^2(\omega)\, d\omega.$$

Using (9) with the substitution $d = j - k$ and Lemma 1(iii), we bound

$$|R_{3n}| \leq \frac{C}{n} \sum_{d=0}^{n-1}\left(d + \sum_{k=1}^{n-1} |a_{k,n}| n^{-1} k (1 - n^{-1}k)\right) \sum_{t \in \mathbb{Z}} |r(t) r(t+d)| \leq CB_n$$

using $\sup_{k,n} |a_{k,n}| \leq C$. Hence, (6) and Theorem 2(i) are established.

For Theorem 2(ii), we use that $|H_n(\omega)| \leq CL_{0,n}(\omega)$, $\omega \in \mathbb{R}$ holds by (5)–(6) in [5] since the weights $|a_{k,n}|$ are bounded and $\sup_n \sum_{k=2}^{n-1} |a_{k,n} - a_{k-1,n}| < \infty$. Using this and that $f^2(\cdot)$ is bounded with $|f^2(0) - f^2(\omega)| \leq C|\omega|$, $\omega \in \mathbb{R}$ by A.1 and that $K_n(\cdot)$ is nonnegative with $\int_{-\pi}^{\pi} K_n(\omega)\, d\omega = 1$, we bound $|\int_{-\pi}^{\pi} K_n(\omega) \times f^2(\omega)\, d\omega - f^2(0)|$ by

$$\int_{|\omega| \leq \varepsilon} + \int_{\varepsilon < |\omega| \leq \pi} K_n(\omega) |f^2(\omega) - f^2(0)|\, d\omega \leq C\varepsilon + CA_n^{-1} L_{0,n}^2(\varepsilon)$$

for any $0 < \varepsilon < \pi$ and $C$ independent of $\varepsilon$. Taking limits as $n \to \infty$ (i.e., $A_n \to \infty$, $L_{0,n}(\varepsilon) \to \varepsilon^{-1}$) and then letting $\varepsilon \downarrow 0$ gives the result.

**Acknowledgments.** The author wishes to thank an Associate Editor and two referees for helpful suggestions that improved the manuscript.

DEPARTMENT OF STATISTICS
IOWA STATE UNIVERSITY
AMES, IOWA 50011
USA
E-MAIL: dnordman@iastate.edu